\documentclass[leqno,12pt]{amsart} 
\usepackage{amssymb}
\usepackage{amssymb,amsmath,amsfonts,amscd}
\usepackage{psfrag,graphicx}
\usepackage{epsfig}
\usepackage[all,dvips,cmtip]{xy}
\newcommand{\IR}{\mathbb{R}}
\newcommand{\IN}{\mathbb{N}}
\newcommand{\IZ}{\mathbb{Z}}

\newcommand{\IC}{\mathbb{C}}

\newcommand{\IF}{\mathbb{F}}
\newcommand{\IP}{\mathbb{P}}

\newcommand{\ii}{\mathrm{i}}

\newcommand{\str}{\star_{\rho}}

\newcommand{\cM}{\mathcal M}
\newcommand{\cO}{\mathcal O}

\newcommand{\cX}{\mathcal X}

\newcommand{\mfp}{\mathfrak{p}}

\newcommand{\al}{\alpha}
\newcommand{\be}{\beta}
\newcommand{\cd}{\cdot}

\newcommand{\e}{\epsilon}
\newcommand{\fr}{\frac }
\newcommand{\ga}{\gamma}
\newcommand{\Ga}{\Gamma}

\newcommand{\lan}{\langle}

\newcommand{\mb}{\mbox}

\newcommand{\nf}{\normalfont}

\newcommand{\op}{\oplus}

\newcommand{\ran}{\rangle}
\newcommand{\ra}{\rightarrow}
\newcommand{\si}{\sigma}
\newcommand{\Si}{\Sigma}

\theoremstyle{plain} 
\newtheorem{theorem}{\indent\sc Theorem}[section]
\newtheorem{lemma}[theorem]{\indent\sc Lemma}
\newtheorem{corollary}[theorem]{\indent\sc Corollary}
\newtheorem{proposition}[theorem]{\indent\sc Proposition}

\theoremstyle{definition} 

\newtheorem{remark}[theorem]{\indent\sc Remark}

\begin{document}

\title[The cohomological crepant resolution conjecture for $\IP(1,3,4,4)$]
{The cohomological crepant resolution conjecture for $\IP(1,3,4,4)$} 

\author{Samuel Boissi\`ere, \'Etienne Mann and Fabio Perroni }

\subjclass[2000]{ 
Primary 14M25; Secondary 14A20.
}

\thanks{ Thanks. This research was partially supported by SNF, 
No 200020-107464/1.
}
\address{
Laboratoire J.A.Dieudonn\'e UMR CNRS 6621 \endgraf
Universit\'e de Nice Sophia-Antipolis \endgraf
Parc Valrose, 06108 Nice \endgraf
France} 
\email{sb@math.unice.fr} 

\address{
D\'epartement de Math\'ematiques \endgraf                                                   
Universit\'e de Montpellier 2   \endgraf                                                  
Place Eug\`ene Bataillon    \endgraf                                                      
F-34 095 Montpellier CEDEX 5  }
\email{emann@math.univ-montp2.fr}

\address{
Institut f\"ur Mathematik \endgraf
Universit\"at Z\"urich \endgraf
Winterthurerstrasse 190 \endgraf
8057 Z\"urich \endgraf
Switzerland
}
\email{fabio.perroni@math.unizh.ch}


\maketitle

\begin{abstract}
We prove the cohomological crepant resolution conjecture
of Ruan for the weighted projective space $\IP(1,3,4,4)$.
To compute the quantum corrected cohomology ring 
we combine the results of Coates-Corti-Iritani-Tseng on $\IP(1,1,1,3)$
and our previous results.
\end{abstract}

\section{Introduction}
The \textit{Cohomological Crepant Resolution Conjecture}, as proposed by Ruan,
predicts the existence of an isomorphism between the Chen-Ruan cohomology ring
of a complex Gorenstein orbifold $ {\cX}$ and the quantum corrected cohomology
ring of any crepant resolution $\rho :Z\rightarrow X$ of the coarse moduli space
$X$ of $ {\cX}$ \cite{Ru}.  The quantum corrected cohomology ring of $Z$ is the
ring obtained from the small quantum cohomology of $Z$ after specialization of
the quantum parameters corresponding to the rational exceptional curves to
$c_1,\dots,c_m$ and the remaining parameters to zero, 
it is denoted by $H^\star_\rho(Z)(c_1,\dots,c_m)$ (see Sec. \ref{sec:qccr}). 
It is an important issue (see \cite{Ru}) to determine the $c_1,\dots,c_m$.
Examples suggest that there can be different choices (see e.g. \cite{BGP}, \cite{CCIT2},
\cite{P}); however, if one assumes the validity of Conjecture 4.1 in \cite{CR},
then Ruan's conjecture follows (Theorem 7.2 in \cite{CR}) 
and the $c_1,\dots, c_m$ acquire a precise meaning (see section 11 in \cite{CR}).

In this paper we consider the weighted projective
space ${\cX}=\IP(1,3,4,4)$. The quantum corrected cohomology ring 
depends on four quantum parameters: $q_1,q_2,q_3$ and $q_4$;
the first three correspond to the components of the exceptional divisor 
over the transverse singularity, while the fourth corresponds to the 
component of the exceptional divisor over the isolated singularity
(see Sec. \ref{sec:qccr}). Our main result, Theorem \ref{th2},
states that, for $(c_1,\dots,c_4)\in \{(\ii,\ii,\ii,1), (-\ii,-\ii,-\ii,1)\}$,
the quantum corrected cohomology ring of $Z$ is isomorphic to 
the Chen-Ruan cohomology ring of $\IP(1,3,4,4)$.
These values of the 
quantum parameters are relevant for the   
cohomological crepant resolution conjecture and its generalization,
Conj. 4.1 in \cite{CR}. In the previous paper \cite{BMP}, 
we proved that $H^\star_{\rho}(Z)(q_1,\dots , q_4)$
becomes isomorphic to $H^\star_{\rm CR}(\IP(1,3,4,4))$ when we set
the quantum parameters to $(\ii,\ii,\ii,0)$ or $(-\ii,-\ii,-\ii,0)$,
which we found strange in regards to the conjecture, where 
the value $c_4=0$ is not considered.
A motivation for this work was that to clarify this point.

To prove the theorem, we compute explicitly a presentation 
of  $H^\star_{\rho}(Z)(c_1,\dots , c_4)$ with quantum parameters 
equal to $(\ii,\ii,\ii,1)$ and $(-\ii,-\ii,-\ii,1)$, then we
give an explicit isomorphism 
$H^\star_\rho(Z;\IC)(c_1,\dots,c_4) \ra H^\star_{\rm CR}(\IP(1,3,4,4))$.
The quantum corrections
coming from the exceptional divisor over the transverse singularity
have been computed in \cite{BMP}.
To compute the quantum corrections coming from the component
of the exceptional divisor over the isolated singularity
we use results from \cite{CCIT}.

\medskip
The paper is organized as follows. In Sec. \ref{sec:bg} we collect 
some background material: we first 
give a presentation of the Chen-Ruan cohomology ring of $\IP(1,3,4,4)$;
then we describe the crepant resolution of $|\IP(1,3,4,4)|$;
finally we write a presentation of $H^\star_{\rho}(Z)(c_1,\dots,c_4)$
with $(c_1,\dots,c_4)=(\ii,\ii,\ii,c_4)$, here $c_4$ is a generic complex number
and we assume the convergence of the power series \eqref{epsilon}.
The fact that \eqref{epsilon} converges in a neighborhood of the origin
to an analytic function that admits analytic continuation in $1\in \IC$
is proved in Sec. \ref{sec:ccit}; more precisely we show that \eqref{epsilon}
is equal to a constant structure of the small quantum cohomology of 
the crepant resolution of $|\IP(1,1,1,3)|$ and then use results from 
\cite{CCIT}. In the last section \ref{sec:main} we prove Theorem \ref{th2}.
 
\section{Background}\label{sec:bg}
The weighted projective space $\IP(1,3,4,4)$
is the quotient stack  $\left[ \IC^4 - \{ 0\}/\IC^\star\right]$,
where $\IC^\star$ acts diagonally with weights $w_0=1,w_1=3,w_2=4$ and $w_3=4$,
it will be denoted by $\cX$.
The coarse moduli space $X:=|\IP(1,3,4,4)|$ is a projective
variety whose singular locus is the disjoint union of
the curve $C := \{[0:0:x_2:x_3] \} \subset X$
and the isolated point ${\rm P}:= [0:1:0:0] \in X$.
Along $C$, $X$ has transverse ${\rm A}_3$ singularities (see \cite{P}),
the point $\rm P$ is a singularity of type $\frac{1}{3}(1,1,1)$, 
according to  Reid's notation \cite{YPG}. 

\subsection{The Chen-Ruan cohomology}
To compute the Chen-Ruan cohomology ring $H^{\star}_{\rm CR}(\cX;\IC)$
we follow \cite{BMPmodel}. The twisted sectors are indexed by the set
${\rm T}:= \left\{ {\rm exp}(2\pi \ii \ga)\, | \, 
\ga \in \left\{ 0, \frac{1}{3}, \frac{2}{3}, \frac{1}{4}, \frac{1}{2},
\frac{3}{4}\right\} \right\}$. For any $g \in {\rm T}$,
$\cX_{(g)}$ is a weighted projective space: set 
${\rm I}(g):=\left\{ i\in \{ 0,1,2,3 \} \, \mid \, g^{w_i} =1  \right\}$,
then $\cX_{(g)}=\IP(w_{{\rm I}(g)})$, where
$(w_{{\rm I}(g)})=(w_i)_{i\in {\rm I}(g)}$. The inertia stack is the disjoint union 
of the twisted sectors:
$$
{\rm I}\cX = \sqcup_{g\in {\rm T}}\IP(w_{{\rm I}(g)}).
$$ 
As a vector space, the Chen-Ruan cohomology is
the cohomology of the inertia stack; the graded structure is obtained
by shifting the degree of the cohomology of any twisted sector by 
twice the corresponding age, ${\rm age}(g)$.
We have
\begin{eqnarray}\label{crcohomology}
H^p_{\rm CR}(\cX;\IC) &=&
 \op_{g\in {\rm T}}H^{p-2{\rm age}(g)}(\IP(w_{{\rm I}(g)});\IC)\\
&=&H^p(\IP(1,3,4,4);\IC) \op H^{p-2}(\IP(3);\IC) \op H^{p-4}(\IP(3);\IC)\op
\nonumber \\
&& H^{p-2}(\IP(4,4);\IC)\op H^{p-2}(\IP(4,4);\IC)\op H^{p-2}(\IP(4,4);\IC).
\nonumber
\end{eqnarray}
A basis of $H^{\star}_{\rm CR}(\cX;\IC)$ is easily obtained in the following
way: 
set 
$$
H, E_1, E_2, E_3, E_4\in H^{2}_{\rm CR}(\cX;\IC) 
$$
be the image of 
$
c_1 \left(\cO_{\cX}(1)\right) \in H^{2}(\cX;\IC),
$
$1\in H^0(\cX_{\left({\rm exp}(\pi\ii /2)\right)};\IC)$,
$1\in H^0(\cX_{\left({\rm exp}(\pi\ii )\right)};\IC)$,
$1\in H^0(\cX_{\left({\rm exp}(\pi\ii 3/2)\right)};\IC)$
and $1\in H^0(\cX_{\left({\rm exp}(2\pi\ii /3)\right)};\IC)$ respectively,
under the inclusion $H^{\star}(\IP(w_{{\rm I}(g)}))\ra H^\star_{\rm CR}(\cX)$
determined by the decomposition \eqref{crcohomology}.
As a commutative $\IC$-algebra, the Chen-Ruan cohomology ring is generated by 
$H,E_1,E_2,E_3,E_4$ with relations (see \cite{BMPmodel}):
\begin{align*}  
&& HE_{4}, \ E_{1}E_{1}-3HE_{2}, \ E_{1}E_{2}-3HE_{3}, \  
E_{1}E_{3}-3H^{2}, \\  &&
E_{2}E_{2}-3H^{2}, \ E_{2}E_{3}-HE_{1}, \ E_{3}E_{3}-HE_{2}, \  
16H^{3}-E_4^{3}, \\ &&
H^2E_1, \, H^2E_2, \, H^2E_3, \, E_1E_4, \, 
E_2E_4, \, E_3E_4. 
\end{align*}  
We see that the following elements form a basis of 
$H^\star_{\rm CR}(\cX;\IC)$ which we fix for the rest 
of the paper:
\begin{equation}\label{basep}
1,\, H,\, E_1,\, E_2,\, E_3,\, E_4,\, H^2,\, HE_1,\, HE_2,\, HE_3,\, E_4^2,\, H^3.
\end{equation}

\begin{remark}\nf
Note that the elements of our basis 
are different from those used in \cite{BMPmodel}
by a combinatorial factor.

Other methods are suitable in order to compute the Chen-Ruan cup 
product of weighted projective spaces, here are a few:
the results in \cite{BCS}
provide a presentation of the Chen-Ruan cohomology ring
for a general toric Deligne-Mumford stack; results from \cite{M} and from \cite{CCLT}.
\end{remark}

\subsection{The crepant resolution}
We study some properties of the crepant resolution of $X:= | \IP(1,3,4,4) |$.
We begin with the following
\begin{proposition}
The variety $X= |\IP(1,3,4,4) |$ has a unique crepant resolution
$\rho :Z \ra X$, up to isomorphism.
\end{proposition}
\begin{proof}
This is a direct consequence of the following facts: 
the $3$-fold singularity  $\frac{1}{3}(1,1,1)$ has a unique
crepant resolution (see e.g. \cite{FN})
and any variety with transverse ${\rm ADE}$ singularities has a unique crepant
resolution (\cite{P}, Prop. 4.2), up to isomorphism.
\end{proof}

\begin{figure}[!ht]  
\begin{center}  
\psfrag{A}{$v_1$} \psfrag{B}{$v_{2}$}
\psfrag{C}{$v_{3}$} \psfrag{D}{$v_{0}$}  
\psfrag{E}{$Q_1$} \psfrag{F}{$Q_2$} \psfrag{G}{$Q_3$} \psfrag{H}{$Q_4$}  
 \includegraphics[width=0.4\linewidth]{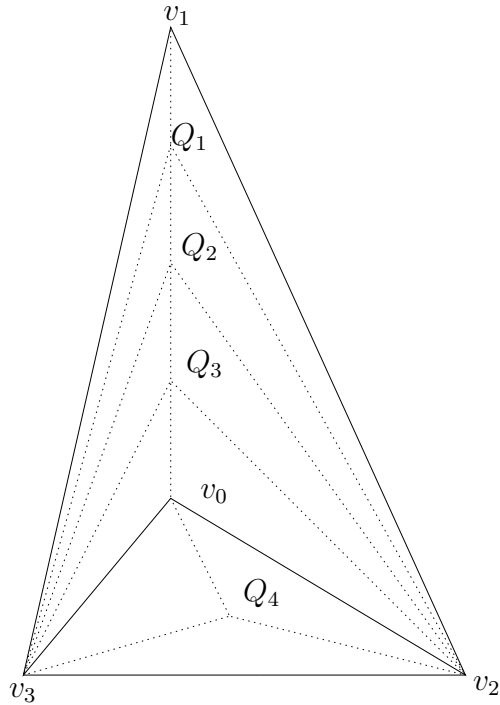}  
\end{center}\caption{Polar polytopes of $|\IP(1,3,4,4)|$ and crepant resolution}\label{fig4}  
\end{figure}

An explicit model for the crepant resolution $\rho :Z \ra X$ can be constructed 
using methods from toric geometry, we follow \cite{F}. 
The toric variety $X$ is associated to the lattice
$\IZ^3$ and the fan $\Si$, where $\Si$ is the fan whose cones are generated by
proper subsets of $\{ v_0:=(-3,-4,-4), v_1:=(1,0,0), v_2=(0,1,0), v_3=(0,0,1)\}\subset \IZ^3$.
The $3$-dimensional cones which correspond to singular open affine subvarieties of $X$
are: $\si_1 := \lan v_0 , v_2 , v_3 \ran$, $\si_2 := \lan v_0 , v_1 , v_3 \ran$, 
$\si_3 := \lan v_0 , v_1 , v_2 \ran$. More precisely, $\si_1$ gives rise to
the isolated singularity, while $\si_2$ and $\si_3$ to the transverse one.
To resolve the isolated singularity, we subdivide $\si_1$ by inserting
the ray generated by 
$$
Q_4:=(-1,-1,-1)= \fr{1}{3}v_0 +\fr{1}{3}v_2 +\fr{1}{3}v_3.
$$
To resolve the transverse singularity, we subdivide $\si_2$ and $\si_3$
by inserting the rays generated by 
\begin{eqnarray*}
Q_1&:=&(0,-1,-1)=\fr{1}{4}v_0 + \fr{3}{4}v_1,\\
Q_2&:=&(-1,-2,-2)=\fr{1}{2}v_0 + \fr{1}{2}v_1, \\
\mb{and} \qquad Q_3&:=&(-2,-3,-3)=\fr{3}{4}v_0 + \fr{1}{4}v_1.
\end{eqnarray*}
The subdivision $\Si'$ is shown in Figure \ref{fig4}.
Then we set 
$Z$ to be the toric variety associated to $\IZ^3$ and $\Si'$,
and $\rho :Z\ra X$ to be the morphism associated to 
the identity $\IZ^3 \ra \IZ^3$. 
\begin{lemma}
The morphism $\rho :Z \ra X$ defined above is a crepant resolution.
\end{lemma}
\begin{proof}
Since any cone of $\Si'$ is generated by a part of a basis of the lattice,
it follows that $Z$ is smooth, \cite{F} Sec. 2.1. 
The crepancy of $\rho$ follows from the existence
of a continuous $\Si$-piecewise linear function  $h':\IR^3 \ra \IR$ such that
$h'(Q_i)=h'(v_j)=-1$, for any $i\in\{1,2,3,4\}$ and $j\in \{0,1,2,3\}$, \cite{F} Sec. 3.4.
\end{proof}
Observe that $\Si'$ is the maximal projective subdivision of the polytope
of $|\IP(1,3,4,4)|$ (see \cite{CK} or \cite{Voisin}).

\subsection{Quantum corrected cohomology ring of $Z$}\label{sec:qccr}
Let us denote by 
$b_j$ (resp. $e_i$) the first Chern class of the holomorphic line bundle 
associated to the torus invariant divisor in $Z$ corresponding to the 
ray generated by $v_j$ (resp. $Q_i$), for any 
$j\in\{ 0,1,2,3\}$  ($i\in\{1,2,3,4\}$ resp.).
Set $h=\fr{1}{12}\left(\sum_{i=0}^3 b_i + \sum_{j=1}^4 e_j\right)$, then
$H^\star(Z;\IC)\cong \IC[h,e_{1},e_{2},e_{3},e_{4}]/I$ where $I$ is the ideal generated by 
(see \cite{BMP} for more details)  
\begin{align*}  
3he_{4}, \   e_{1}e_{3}, \ e_{1}e_{4}, \ e_{2}e_{4}, \ e_{3}e_{4}, \\  
e_{1}^{2} - 10 he_{1}-4he_{2}-2he_{3}+24h^{2},\\  
e_{1}e_{2} + 3 he_{1}+2 he_{2}+he_{3}-12h^{2},\\  
e_{2}^{2}-6he_{1}-12he_{2}-2he_{3}+24h^{2},\\  
e_{2}e_{3}+3he_{1}+6he_{2}+he_{3}-12h^{2},\\  
e_{3}^{2}-6he_{1}-12he_{2}-14he_{3}+24h^{2},\\  
16h^{2}e_{1},\ 16h^{2}e_{2}, \ 16h^{2}e_{3}, 16h^{3}-\frac{1}{27}e_{4}^{3}.  
\end{align*}  
Let us fix the following basis of  $H^\star(Z;\IC)$:
\begin{eqnarray}\label{basez}
\psi_0=1,\, \psi_1=h,\, \psi_2=e_1,\, \psi_3=e_2,\, \psi_4=e_3,\, 
\psi_5=e_4,\\
\psi_6=h^2,\, \psi_7=he_1,\, \psi_8=he_2,\, \psi_9=he_3,\, \psi_{10}=e_4^2,\, \psi_{11}=h^3,\nonumber
\end{eqnarray}
with dual basis:
\begin{eqnarray*}
\psi^0=48h^3,\, \psi^1=48h^2,\, \psi^2=-3he_1 -2he_2 -he_3,\, 
\psi^3=-2he_1-4he_2-2he_3,\\ \psi^4=-he_1-2he_2-3he_3,\, 
\psi^5=\fr{1}{9}e_4^2,\,  \psi^6=48h,\, \psi^7=-3e_1-2e_2-e_3,\\ 
\psi^8=-2e_1-4e_2-2e_3,\, \psi^9=-e_1-2e_2-3e_3,\, \psi^{10}=\fr{1}{9}e_4,\, \psi^{11}=48\cd 1.
\end{eqnarray*}

Let ${\rm M}_{\rho}(Z)\subset {\rm A}_1(Z;\IZ)$ be the cone of effective $1$-cycles
in $Z$ which are contracted by $\rho$. It is freely generated by
$\Ga_i :={\rm PD}(4he_i)$ for $i\in \{1,2,3\}$
and $\Ga_4:={\rm PD}\left( -\fr{1}{3}e_4^2 \right)$, where $\rm PD$ means 
Poincar\'e dual (see e.g. \cite{BMP}). Let $q_1,q_2,q_3,q_4$ be formal variables and let
$\Lambda := \IC [\![ q_1,\dots,q_4 ]\!]$ be the ring of 
formal power series in $q_1,...,q_4$. We have an associative product on
$H^\star(Z;\Lambda)$ defined as:
\begin{equation}\label{starrho}
\psi_i \star_{\rho} \psi_j := \sum_{\Ga \in {\rm M}_{\rho}(Z)} 
\sum_{\ell=0}^{11}\lan \psi_i,\psi_j, \psi^\ell\ran^Z_{0,3,\Ga} \psi_\ell q^\Ga,
\end{equation}
where $\lan \dots \ran^Z_{0,3,\Ga}$ is the Gromov-Witten invariant of $Z$ of genus zero, 
three marked points, homology class $\Ga$ and $q^\Ga:=q_1^{d_1}\cd \cd \cd q_4^{d_4}$
for $\Ga = d_1\Ga_1 +...+ d_4 \Ga_4 \in {\rm M}_{\rho}(Z)$. 

In \cite{BMP} we computed explicitly the product \eqref{starrho} 
whenever $\psi_i \not= e_4$ or $\psi_j \not= e_4$; as a result it follows that,
in these cases, the power series involved in \eqref{starrho} converge in a neighborhood 
of the origin. We will see in Section \ref{sec:ccit} that also  
$e_4\star_{\rho} e_4$ is convergent in a neighborhood of the origin.
Therefore the expression \eqref{starrho} defines a family of $\IC$-algebras over the vector space 
$H^\star(Z;\IC)$ whose structure constants are analytic functions defined
in some region of the complex space $\IC^4$. For any $c_1,\dots,c_4$
in this region, the algebra obtained by setting $q_i=c_i$ 
is the \textit{quantum corrected cohomology ring}
of $Z$ with quantum parameters specialized to $c_1,\dots,c_4$ and is denoted by
$H^\star(Z;\IC)(c_1,\dots,c_4)$ \cite{Ru}.

In particular we have the following presentation for
$H^\star_\rho(Z)(\ii,\ii,\ii,c_4)$ (see \cite{BMP} for more details):
\begin{align}\label{qcproduct}
&\al_1 \str \al_2 = \al_1 \cup \al_2, \quad \mb{if} \quad 
{\rm deg}(\al_1) \not= 2,\nonumber\\
&h\str \al = h\cup \al, \quad \mb{for any} \, \al,\nonumber \\
&e_{1}\str e_{1}=  
-24h^{2}+(-2+6\ii)he_{1}-4he_{2}+(-2-2\ii)he_{3},\nonumber\\  
& e_{1}\str e_{2}=12h^{2}+(-1-4\ii)he_{1}+(2-4\ii)he_{2}+he_{3},\nonumber\\  
& e_{1}\str e_{3}=  
-2\ii he_{1}-2\ii he_{3}, \\  
& e_{2}\str e_{2}=  
-24h^{2}+(2+2\ii)he_{1}+8\ii he_{2}+(-2+2\ii)he_{3},\nonumber\\  
& e_{2}\str e_{3}=  
12h^{2}-he_{1}+(-2-4\ii)he_{2}+(1-4\ii)he_{3},\nonumber\\  
& e_{3}\str e_{3}=  
-24h^{2}+(2-2\ii)he_{1}+4he_{2}+(2+6\ii)he_{3}, \nonumber \\
& e_4 \str e_4 = \e(c_4) e_4^2; \nonumber
\end{align}  
where $\cup$ is the usual cup product, 
and 
\begin{eqnarray}\label{epsilon}
\e (q) &:=& 1+\fr{1}{9}\left(\int_{\Ga_4}e_4\right)^3\sum_{a=1}^\infty a^3
{\rm deg}\left[\overline{\cM}_{0,0}(Z,a\Ga_4)\right]^{\rm vir}{q}^a \\
&=& 1-3\sum_{a=1}^\infty a^3
{\rm deg}\left[\overline{\cM}_{0,0}(Z,a\Ga_4)\right]^{\rm vir}{q}^a. \nonumber
\end{eqnarray}

A similar presentation holds for $H^\star(Z)(-\ii,-\ii,-\ii,c_4)$, also 
in this case $ e_4 \str e_4 = \e(c_4) e_4^2$ with $\e(q)$ defined in 
\eqref{epsilon}.

\section{Relations with the quantum  cohomology of $\IF_3$}\label{sec:ccit}
To state our main result (Theorem \ref{th2}), we need to know that the series \eqref{epsilon}
converges in a neighborhood of the origin to an analytic function that admits analytic 
continuation at the point $1\in \IC$. To this aim, we show in this section
that \eqref{epsilon} appears as a structure constant of the small quantum cohomology 
of $\IF_3$, the crepant resolution of $|\IP(1,1,1,3)|$,
so that we can use the results of \cite{CCIT} and \cite{CR}.
For reader's convenience, we stick to the notation of \cite{CCIT}.

\medskip
Let us consider the coarse moduli space of the 
weighted projective space $\IP(1,1,1,3)$,
it has an isolated singularity of type $\fr{1}{3}(1,1,1)$.
Let $\IF_3:=\IP(\cO_{\IP^2}(-3)\op \cO_{\IP^2})$, then there 
exists a morphism $\chi:\IF_3 \ra |\IP(1,1,1,3)|$
which is a crepant resolution. 
Following the notation of \cite{CCIT},
let $p_1 \in H^2(\IF_3;\IC)$ be the class Poincar\'e dual to the preimage 
in $\IF_3$ of a hyperplane in $\IP^2$;
let $p_2 \in H^2(\IF_3;\IC)$ the class Poincar\'e dual to the infinity section.
The cohomology ring of $\IF_3$ has the following presentation 
$$
H^\star(\IF_3;\IC)=\IC [p_1,p_2]/\lan p_1^3, p_2^2 -3p_1p_2\ran.
$$
Let us fix the following basis of $H^\star(\IF_3;\IC)$ (as in \cite{CCIT})
\begin{equation*}
\phi_0=1, \, \phi_1=\fr{p_2}{3}, \, \phi_2=\fr{p_1p_2}{3}, \, \phi_3=\fr{p_2-3p_1}{3}, \,
\phi_4=-\fr{p_1(p_2-3p_1)}{3}, \, \phi_5=\fr{p_1^2p_2}{3},  
\end{equation*}
then its dual is 
$$
\phi^0=p_1^2p_2, \, \phi^1=p_1p_2, \, \phi^2=p_2, \, \phi^3= -p_1(p_2-3p_1),
\, \phi^4= p_2-3p_1, \, \phi^5=3.
$$
The Poincar\'e dual
to the exceptional divisor is $\mfp := p_2 -3p_1$ (in \cite{CCIT}
$\mfp$ is denoted by $\mfp_1$).  
We want to compute the product $\mfp \circ_{\hat q}\mfp$ in the small
quantum cohomology of $\IF_3$ with $(\hat{q}_1,\hat{q}_2)=(1,0)$ (see \eqref{decq}). 

By definition, we have (see \cite{CCIT}, sec. 2.4)
\begin{equation}\label{pop}
\mfp \circ_{\hat q}\mfp =\sum_{\ga \in {\rm M}(\IF_3)}\sum_{\ell =0}^5
\lan \mfp,\mfp,\phi_\ell \ran^{\IF_3}_{0,3,\ga} \hat{q}^{\ga}\phi^\ell,
\end{equation}
where ${\rm M}(\IF_3)\subset H_2(\IF_3;\IZ)$ is the cone of the classes of
effective curves, for $\ga \in {\rm M}(\IF_3)$ 
$\hat{q}^\ga : \fr{H^2(\IF_3;\IC)}{2\pi \ii H^2(\IF_3;\IZ)}\ra \IC^\star$
is the function $[\tau] \mapsto {\rm exp}\left( \int_\ga\tau\right)$,
$\lan \dots \ran^{\IF_3}_{0,3,\ga}$ is the Gromov-Witten invariant of $\IF_3$
of genus zero, three marked points and class $\ga$.

Since $p_1,\, p_2$ form a basis of $H^2(\IF_3;\IC)$,  any 
$[\tau]$ can be written as $[\tau]=[\tau^1p_1 +\tau^2p_2]$, therefore
\begin{equation}\label{decq}
\hat{q}^\ga([\tau])= \hat{q}_1^{\int_\ga p_1}\hat{q}_2^{\int_\ga p_2},
\end{equation}
where $\hat{q}_i:={\rm exp}(\tau^i)$, $i\in\{1,2\}$.
Let us consider the classes $\ga_1:={\rm PD}(-\fr{1}{3}\mfp^2=p_1(p_2-3p_1))$
 and $\ga_2:={\rm PD}(p_1^2)$, they form a basis for ${\rm M}(\IF_3)$
such that $\int_{\ga_i}p_j=\delta_{ij}$, $i,j \in \{1,2\}$.
It follows that the product $\mfp \circ_{\hat q}\mfp$
restricted to $\hat{q}=(\hat{q}_1,0)$ is given by:
\begin{equation*}
\mfp \circ_{(\hat{q}_1,0)}\mfp =\sum_{a\geq 0}\sum_{\ell =0}^5
\lan \mfp,\mfp,\phi_\ell \ran^{\IF_3}_{0,3,a\ga_1} \hat{q}_1^a\phi^\ell.
\end{equation*}
Using the degree axiom, the divisor axiom and $\int_{\ga_1}p_2=0$, we deduce that 
 \begin{eqnarray*}
\mfp \circ_{(\hat{q}_1,0)}\mfp =&&\sum_{a =0}^\infty
\lan \mfp,\mfp,\fr{1}{3}\mfp \ran^{\IF_3}_{0,3,a\ga_1}\cd (-p_1(p_2-3p_1))\hat{q}_1^a.
\end{eqnarray*}
Finally, since $\mfp^2:=\mfp \cup \mfp = -3p_1(p_2-3p_1)$, using again the divisor axiom
we rewrite the previous expression as
\begin{equation}\label{ultimo}
\mfp \circ_{(\hat{q}_1,0)}\mfp = \left( 1 -3 \sum_{a=1}^\infty a^3 
{\rm deg}\left[\overline{\cM}_{0,0}(\IF_3,a\ga_1)\right]^{\rm vir} \hat{q}_1^a \right)\mfp^2.
\end{equation}

Then we have:
\begin{lemma}\label{e=}
The power series in \eqref{ultimo} is equal to the power series $\e$
(see Formula \eqref{epsilon}). 
\end{lemma} 
\begin{proof}
From the uniqueness of the crepant resolution of the singularity
$\fr{1}{3}(1,1,1)$, there exists an isomorphism
between a neighborhood in $Z$ of the component of the 
exceptional divisor over $[0:1:0:0]\in |\IP(1,3,4,4)|$
and a neighborhood in $\IF_3$ of the exceptional divisor of $\chi$. 
The isomorphism between these neighborhoods
induces isomorphisms of stacks (see e.g. \cite{P},
Lemma 7.1):
\begin{equation}\label{ms}
\overline{\cM}_{0,0}(Z,a\Ga_4) \cong  \overline{\cM}_{0,0}(\IF_3,a\ga_1) \qquad 
\forall a\in \IN-\{0\}. 
\end{equation}
Since the virtual fundamental classes
of the above stacks  depend  only on neighborhoods
of the exceptional divisors, it follows that
$$
{\rm deg}\left[\overline{\cM}_{0,0}(Z,a\Ga_4)\right]^{\rm vir}=
{\rm deg}\left[\overline{\cM}_{0,0}(\IF_3,a\ga_1)\right]^{\rm vir}  \qquad \forall a\in \IN-\{0\}.
$$
Hence the result follows.
\end{proof} 
As a consequence we have:
\begin{corollary}\label{e1}
The series \eqref{epsilon} converges in a neighborhood of the origin
to an analytic function $f(q)$ that has analytic continuation in $1\in \IC$. 
Moreover the value of $f(q)$ at $q=1$ is
$$
f(1)=\fr{2\pi \be_1}{9\be_2^2}=\frac{(2\pi)^6}{27\Ga\left(\fr{1}{3}\right)^9},
$$
where $\be_i=\fr{2\pi}{9\Ga\left(\fr{i}{3}\right)^3}$, $i\in \{1,2\}$.
\end{corollary}
\begin{proof}
It follows from \cite{CCIT} (see also \cite{CR} Theorem 7.2)
that the statement is true for \eqref{ultimo} and Lemma \ref{e=}
implies that it is also true for \eqref{epsilon}.

To compute $f(1)$, consider the map $\Theta(q):H^\star_{\rm CR}(\IP(1,1,1,3))\ra H^\star(\IF_3)$
defined in \cite{CCIT} (pag. 56) and set $q=0$.
We consider $H^\star(\IF_3)$ with the algebra structure coming
from the quantum cohomology when we set the quantum parameter $\hat{q}$ 
at $(1,0)$ (see \eqref{decq}).
$\Theta(0)$ is a morphism of  algebras, 
so we have \cite{CCIT}\footnote{The careful reader will notice that, in \cite{CCIT},     
$\IP(1,1,1,3)$ is defined as the quotient  $[\IC^4- \{ 0 \}/\IC^\star]$ where the action    
has negative weights.}:
$$
\fr{2\pi}{3}\be_1 \mfp^2 =\Theta(0)(\mathbf{1}_{\fr{1}{3}})
=\Theta(0)(\mathbf{1}_{\fr{2}{3}})\circ_{(1,0)}\Theta(0)(\mathbf{1}_{\fr{2}{3}})=
3\be_2^2 \mfp \circ_{(1,0)} \mfp, 
$$
therefore 
$$
\mfp \circ_{(1,0)}\mfp = f(1)\mfp^2=\fr{2\pi  \be_1}{9\be_2^2} \mfp^2 \quad
\mb{and} \quad f(1)=\fr{2\pi \be_1}{9\be_2^2}.
$$
The second equality follows from standard identities of the $\Ga$-function.
\end{proof}

\section{The main result}\label{sec:main}
From the results in \cite{BMP} and Cor. \ref{e1} it follows that 
the power series in \eqref{starrho} converge in a neighborhood 
of the origin of $\IC^4$ to analytic functions that admit 
analytic continuations in $(\ii,\ii,\ii,1)$ and in $(-\ii,-\ii,-\ii,1)$,
moreover the following result holds.
\begin{theorem} \label{th2} 
  For  $(c_1,c_2,c_3,c_4) \in \{ (\ii,\ii,\ii,1), (-\ii,-\ii,-\ii,1) \}$ 
  there is a ring isomorphism 
\begin{equation*} 
H^{\star}_{\rho}(Z;\IC)(c_1,c_2,c_3,c_4) \stackrel{\cong}{\longrightarrow}
 H^{\star}_{\rm CR}(\IP(1,3,4,4);\IC) 
\end{equation*} 
which is an isometry with respect to the Poincar\'e pairings.
\end{theorem} 
\begin{proof} Let us first consider the case where $(c_1,\dots,c_4) = (\ii,\ii,\ii,1)$.
Let  
\begin{equation}\label{mappa}
\Xi : H^\star_\rho(Z;\IC)(\ii,\ii,\ii,1) \longrightarrow H^\star_{\rm CR}(\IP(1,3,4,4);\IC)
\end{equation}
be the linear map defined by the following matrix with respect to the basis \eqref{basep}
and \eqref{basez},
\begin{equation*} \left(
\begin{array}{cccccccccccc}
1 & 0 & 0            & 0       &        0     &   0 & 0 & 0            & 0     & 0            & 0   & 0  \\
0 & 1 & 0            &     0   &        0     &   0 & 0 & 0            & 0     & 0            & 0   & 0  \\
0 & 0 & -\sqrt{2}    & -2\ii   & \sqrt{2}     &   0 & 0 & 0            & 0     & 0            & 0   & 0\\
0 & 0 & -\ii\sqrt{2} & 2\ii    & -\ii\sqrt{2} &   0 & 0 & 0            & 0     & 0            & 0   & 0 \\
0 & 0 & \sqrt{2}     & -2\ii   & -\sqrt{2}    & 0   & 0 & 0            & 0     & 0            & 0   & 0\\
0 & 0 & 0            & 0       & 0            & \al & 0 & 0            & 0     & 0            & 0   & 0\\
0 & 0 & 0            &       0 & 0            &   0 & 1 & 0            & 0     & 0            & 0   & 0\\
0 & 0 & 0            & 0       & 0            & 0   & 0 & -\sqrt{2}    & -2\ii & \sqrt{2}     & 0   & 0\\
0 & 0 & 0            & 0       & 0            & 0   & 0 & -\ii\sqrt{2} & 2\ii  & -\ii\sqrt{2} & 0   & 0 \\
0 & 0 & 0            & 0       & 0            & 0   & 0 & \sqrt{2}     & -2\ii & -\sqrt{2}    & 0   & 0\\
0 & 0 & 0            & 0       & 0            & 0   & 0 & 0            & 0     & 0            & \be & 0\\
0 & 0 & 0            & 0       & 0            & 0   & 0 & 0            & 0     & 0            & 0   & 1
\end{array} \right)
\end{equation*}

\medskip
\noindent where $\al$ and $\be$ are complex numbers to be determined. 
Note that $\Xi$ coincides with the map (4.3) in \cite{BMP}
except possibly for the image of $e_4$ and ${e_4}^2$.

Clearly $\Xi$ is an isomorphism of vector spaces if and only if 
$\al \cd \be \not=0$. We have that
$\Xi(e_i\str e_j)= \Xi(e_i)\cup_{\rm CR}\Xi(e_j)$
for $i,j \in \{1,2,3\}$ \cite{BMP}, therefore it remains to find $\al$
and $\be$ such that 
$\al \cd \be \not=0$, $\Xi(e_4\str e_4)= \Xi(e_4)\cup_{\rm CR}\Xi(e_4)$
and $\Xi(e_4\str e_4^2)= \Xi(e_4)\cup_{\rm CR}\Xi(e_4^2)$.
This is equivalent to the equations:
\begin{equation}\label{gamma}
\al^3 = 27 f(1) \quad \rm{and}\quad \al\cd \be = 27.
\end{equation}

The existence of $\al$ and $\be$ verifying \eqref{gamma} follows from Corollary \ref{e1}.
The fact that $\Xi$ is an isometry is then a direct consequence of \cite{BMP}
and \eqref{gamma}. 

The case where $(c_1,\dots, c_4)=(-\ii,-\ii,-\ii,1)$ is similar
to the previous one, the linear map $H^\star_\rho(Z;\IC)(-\ii,-\ii,-\ii,1) \ra H^\star_{\rm CR}(\IP(1,3,4,4);\IC)$
coincides with (4.4) in \cite{BMP} on $1,e_1,e_2,e_3,h^2,he_1,he_2,he_3,h^3$
and coincides with \eqref{mappa} on $e_4$ and $e_4^2$.
Also in this case the fact that it is a ring isomorphism and that preserves the metrics
follows from \cite{BMP} and \eqref{gamma}. This concludes the proof. 
\end{proof}

\bibliographystyle{amsalpha} 
\bibliography{Biblio.bib} 
 
\end{document}